\documentclass{ifacconf}

\usepackage{graphicx}      
\usepackage{natbib}        
\usepackage{amsmath}
\usepackage{amssymb} 
\begin{document}
\begin{frontmatter}

\title{Revisiting Diffusive Representations for Enhanced Numerical Approximation of Fractional Integrals\thanksref{footnoteinfo}} 
\thanks[footnoteinfo]{This work has been supported by 
	the German Federal Ministry for Education and Research (BMBF)
	under Grant No.\ 05M22WHA.}

\author[First]{Renu Chaudhary} 
\author[First]{Kai Diethelm}

\address[First]{FANG, Technical University of Applied Sciences Würzburg-Schweinfurt, Ignaz-Schön-Str. 11, 97421 Schweinfurt, Germany (Email: \{renu.chaudhary, kai.diethelm\}@thws.de)}

\begin{abstract}                
This study reexamines diffusive representations for fractional integrals with the goal of pioneering new variants of such representations. These variants aim to offer highly efficient numerical algorithms for the approximate computation of fractional integrals. The approach seamlessly aligns with established techniques used in addressing problems involving integer-order operators, contributing to a unified framework for numerical solutions.
\end{abstract}


\end{frontmatter}
\section{Introduction and Preliminaries}
The realm of fractional calculus, with its roots dating back to the late 17th century, has gained significant prominence in various scientific and engineering disciplines in recent decades. Fractional calculus offers a unique and versatile framework for describing complex phenomena that extend beyond the scope of classical calculus. Among its key components, fractional integrals play a pivotal role in understanding and modeling various physical, mathematical, and engineering systems, where memory and non-local behaviors are prevalent.

In the pursuit of harnessing the power of fractional integrals for solving practical problems, one encounters a fundamental challenge that is their numerical evaluation. Unlike their classical integer-order counterparts, fractional integrals exhibit a non-local character that complicates the computation process, making it significantly more resource-intensive in terms of both runtime efficiency and memory requirements. This challenge becomes particularly pronounced when dealing with complex systems or the need to evaluate fractional integrals at multiple points.

However, researchers and practitioners have sought innovative solutions to overcome these computational hurdles and unlock the potential of fractional integrals in various applications. One promising avenue in this endeavor is the diffusive representation of fractional integrals. This representation introduces a novel perspective on how fractional integrals can be expressed and evaluated, offering a dynamic and sustainable solution to address the computational complexities inherent in these non-local operators. For a more comprehensive understanding of the diffusive representation of fractional operators, we refer to \cite{montseny1998diffusive}, \cite{diethelm2022new}, \cite{diethelm2021new} and references therein. 

In this paper, we embark on an exploration of an innovative adaptation of the diffusive representation method for the Riemann-Liouville integral of order $\alpha\in (0,1)$ pertaining to functions $f\in C[a,b]$ (where $a$ and $b$ are real numbers satisfying $a<b$) commencing from point $a$, denoted as
\begin{equation}\label{RLIntegral}
I^{\alpha}_af(t)=\frac{1}{\Gamma(\alpha)}\int_a^t(t-\tau)^{\alpha-1}f(\tau) \, \mathrm d\tau.
\end{equation}

The diffusive representation of (\ref{RLIntegral}) is structured as an integral that incorporates an auxiliary bivariate function. This integral, notable for its evaluation within a fixed range with respect to a single variable, facilitates efficient numerical computation through the utilization of a predefined quadrature formula, thereby reducing the arithmetic complexity. Moreover, the auxiliary function can be ascertained through the solution of a first-order differential equation, thereby eliminating the necessity to consider memory effects during the computational procedure.

Our aim is to improve existing diffusive representations, offering better algorithms for approximating fractional integrals numerically. In doing so, we provide researchers and practitioners with effective tools for handling the complex calculations involved in fractional integrals. These upgraded methods make it easier to solve difficult problems with greater precision and efficiency, ultimately making fractional calculus more practical and useful in various fields.

We first provide the fundamental definitions and results as outlined in \cite{diethelm2023diffusive}.

\begin{defn}
A function $\psi: \Omega \rightarrow (0,\infty)$ is called an \emph{admissible transformation} if it satisfies the following criteria:
\begin{description}
  \item[(i)] The domain $\Omega$ represents a non-empty open interval.
  \item[(ii)] $\psi$ is strictly monotonically increasing.
  \item[(iii)] $\psi\in C^1(\Omega)$.
  \item[(iv)] $\lim_{r\rightarrow \inf \Omega}\psi(r)=0$ and $\lim_{r\rightarrow \sup \Omega}\psi(r)=+\infty$.
\end{description}
\end{defn}

\begin{thm}
For a given admissible transformation  $\psi: \Omega \rightarrow (0,\infty)$, the Riemann-Liouville integral of order $\alpha>0$, $\alpha\in \mathbb{N}$ for a function $f\in C[a,b]$, with $a, b \in \mathbb R$ and $a<b$, can be expressed as
\begin{equation}\label{RLDiffusiveIntegral}
I^{\alpha}_af(t)=\int_{\Omega}\phi(t,r) \, \mathrm dr
\end{equation}
with
\begin{equation}
\begin{array}{ll}
\phi(t,r)=c_{\alpha}&\psi^{\prime}(r)(\psi(r))^{n-\alpha-1}\\
&\times\int_a^t(t-\tau)^{n-\alpha-1}\mathrm e^{-(t-\tau)\psi(r)}f(\tau) \, \mathrm d\tau,\label{phi}
\end{array}
\end{equation}
where $n=\lceil\alpha\rceil$ and the constant $c_{\alpha}$ is given by 
\begin{equation*}
c_{\alpha}=\frac{\sin\pi \alpha}{\pi}\cdot \prod_{l=1}^{n-1}\frac{1}{l-\alpha}.
\end{equation*}
Moreover, when considering a fixed value of $r \in \Omega$, the function $\phi(\cdot, r)$ is distinguished as the unique solution to the $n$-th order differential equation
\begin{equation}\label{DE}
\begin{array}{ll}
\sum_{k=0}^n \binom{n}{k} &(\psi(r))^{n-k}\frac{\partial^k \phi(t,r)}{\partial t^k}\\
&=c_{\alpha}\psi^{\prime}(r)(\psi(r))^{n-\alpha-1}(n-1)!f(t)
\end{array}
\end{equation}
for $t \in [a,b]$, subject to the initial conditions
\begin{equation}
\frac{\partial^k \phi(a,r)}{\partial t^k}=0.
\end{equation}
\end{thm}

\begin{defn}
	The representation \eqref{RLDiffusiveIntegral} of the Riemann-Liouville integral is called a \emph{diffusive representation} of $I_a^\alpha f$.
\end{defn}

These definitions serve as fundamental concepts in the understanding and application of the admissible transformation and diffusive representation in the context of fractional integrals. 

In the existing literature, various methodologies have been developed to express the fractional integral of a given function. For instance, \cite{RebeccaLi} computed the fractional integral by employing an integral representation of the convolution kernel. Subsequently, she devised an efficient quadrature method for this integral representation and incorporated it into a fast time-stepping technique. In another study by \cite{beylkin}, the authors applied Poisson's summation formula to discretize the integral representation of a power function, enabling them to approximate it using exponential functions. This approach allowed for more efficient computation.

Furthermore, in the research presented by \cite{mclean2018exponential}, the idea proposed by \cite{beylkin} was extended and further developed by introducing an alternative integral representation of the power function. This alternative representation provided additional flexibility and accuracy in fractional integral computation. In a different study by \cite{baffet}, the Gauss-Jacobi quadrature rule was employed to approximate the kernel of the fractional integral through a linear combination of exponentials. This technique offered an effective means of numerical approximation for the fractional integral.

In this paper, our focus is to develop a novel diffusive representation for the integral presented in (\ref{RLDiffusiveIntegral}). This new representation provides a pragmatic pathway to design exceptionally efficient numerical algorithms for approximating the integral (\ref{RLDiffusiveIntegral}) with reduced memory usage and computational complexity. This approach aligns with the conventional techniques commonly employed in solving problems involving integer-order operators. To accomplish this, we have examined a unique form of the admissible transformation function, denoted as $\psi(\cdot)$, specifically tailored to make the associated auxiliary bivariate function $\phi(t,\cdot)$ display exponential decay characteristics as it approaches the boundaries of the domain $\Omega$. This characteristic bears substantial significance in the context of approximation theory, as it allows us to leverage both the truncated trapezoidal scheme and the Gauss-Laguerre numerical quadrature scheme for the integration of the function $\phi(t, \cdot)$. Furthermore, it is essential to possess an understanding of the integrand $\phi(t,\cdot)$, and this understanding has been achieved by numerically solving the initial value problem indicated in (\ref{DE}) through the backward Euler's method and trapezoidal method.

\section{Variants of Diffusive Representation Enabling Efficient Numerical Treatment}

In this section, we discuss a promising variant of the diffusive representation to the integral represented in (\ref{RLDiffusiveIntegral}) for the fractional order $\alpha\in (0,1)$. We illustrate that within this novel representation, the integrand $\phi(t,r)$ demonstrates an exponential decay as the variable $r$ approaches the boundary points of the integration domain. This property, combined with the smoothness property, facilitates notably more efficient numerical integration.

\subsection{The Fundamental Approach}

\begin{thm}\label{mainthm2}
Let $\alpha \in (0,1)$ and
consider a function $f\in C[a,b]$ with some $a, b \in \mathbb R$ such that $a<b$. Define the function $\phi(t,r)$ as
\begin{equation}\label{ExDiffusiveIntegralpart2}
\begin{split}
 \phi(t,r)= {} & c_{\alpha}(1+\mathrm e^{-r})\exp[(1-\alpha)(r-\mathrm e^{-r})]\\
& {} \times\int_a^t\exp[-(t-\tau)\exp(r-\mathrm e^{-r})]f(\tau) \, \mathrm d\tau
\end{split}
\end{equation}
for all $r\in \mathbb{R}$ and $t\in [a,b]$. The ensuing properties can then be observed:
\begin{description}
  \item[(i)] The function $\phi(\cdot, r)$ is the unique solution to the initial value problem for a first-order differential equation on the interval $[a, b]$
\begin{eqnarray}\label{DE2}
\begin{split}
\frac{\partial \phi(t,r)}{\partial t}&=-\exp(r-\mathrm e^{-r}) \phi(t,r) \\
&\qquad {} + c_{\alpha}(1+\mathrm e^{-r})[\exp(r-\mathrm e^{-r})]^{1-\alpha}f(t),\\
\phi(a,r)&=0.
\end{split}
\end{eqnarray}
  \item[(ii)] For any $t\in [a,b]$,
\begin{equation}\label{ExDiffusiveIntegral11}
I^{\alpha}_af(t)=\int_{-\infty}^{\infty}\phi(t,r) \, \mathrm dr.
\end{equation}
  \item[(iii)] For any $t\in [a,b]$, the integrand $\phi(t,\cdot)\in C^{\infty}(\mathbb{R})$.
  \item[(iv)] For any $t \in [a, b]$, there exists a constant $C > 0$ such that
\begin{equation}\label{asymp 3}
 |\phi(t,r)|\leq C (1+\mathrm e^{-r})\mathrm e^{-\alpha (r-\mathrm e^{-r})} \quad \text{as}\,\,\, r \rightarrow \infty
\end{equation}
and 
\begin{equation}\label{asymp 4}
 |\phi(t,r)|\leq C(1+\mathrm e^{-r})\mathrm e^{(1-\alpha) (r-\mathrm e^{-r})} \quad \text{as}\,\,r\rightarrow -\infty.
\end{equation}
\end{description}
\end{thm}

\begin{pf} 
Noting that our assumption $\alpha \in (0,1)$ implies that we have to choose $n=1$ in (\ref{phi}),
the proof of the Theorem can be directly deduced from the results presented in \cite{diethelm2023diffusive}
upon choosing $\psi(r)=\exp(r-\mathrm e^r)$ for $r \in \Omega=(-\infty,\infty)$. 
\qed
\end{pf}

Our decision to use this special transformation $\psi$ has been motivated by observations described by \cite{mclean2018exponential} who noticed some very favourable properties of this transformation in a related context.

To further analyze the features of this transformation, we use the following notation: Given two functions $f$ and $g$ defined over $\mathbb R$ we write that $f(x) \sim g(x)$ for $x \to \infty$ if there exist two positive numbers $A$ and $B$ such that 
\[
	A \le  \frac{f(x)}{g(x)}  \le B
\]
whenever $x$ is sufficiently large (and analogously for $x \to -\infty$). Effectively this means that, as $x$ tends to the indicated limit, the asymptotic behaviour of $f(x)$ is identical to the behaviour of $g(x)$; in particular, both functions decay or increase with identical rates.

\begin{rem}
	\label{rem:freud-erdos}
A close inspection of the decay properties of $|\phi(t, r)|$ for fixed $t \in [a,b]$ and $r \to \pm \infty$ as described in eqs.~\eqref{asymp 3} and
\eqref{asymp 4} reveals the following properties:
\begin{enumerate}
\item[(a)] For $r \rightarrow \infty$, the upper bound for $|\phi(t,r)|$ given in (\ref{asymp 3}) is
	\[
		C (1 + \mathrm e^{-r}) \mathrm e^{-\alpha (r - \mathrm e^{-r})} \sim \mathrm e^{-\alpha r}
	\] 
	because $1 + \mathrm e^{-r} \sim 1$ and $r - \mathrm e^{-r} \sim r$. 
	Therefore---since this is only an upper bound---$|\phi(t, r)|$ decays at least as fast as $\mathrm e^{-\alpha r}$. Indeed, a careful look at the details of the proof of this estimate \citep[proof of Theorem 5]{diethelm2023diffusive} shows that a faster decay rate than $\mathrm e^{-\alpha r}$ is possible but only for rare and exceptional choices of the function $f$.
\item[(b)] However, as  $r \rightarrow -\infty$, we observe a significantly different behaviour. Specifically, in this case we have for the terms appearing on the right-hand side of \eqref{asymp 4} that $1 + \mathrm e^{-r} \sim \mathrm e^{-r}$ and $r-\mathrm e^{-r} \sim - \mathrm e^{-r}$, and therefore the upper bound of \eqref{asymp 4} behaves as
\begin{equation}
	\label{3}
	C(1+\mathrm e^{-r})\mathrm e^{(1-\alpha) (r-\mathrm e^{-r})} 
	\sim \mathrm e^{-r -(1-\alpha) \mathrm e^{-r}}
	\sim \mathrm e^{-(1-\alpha) \mathrm e^{-r}}.
\end{equation}
\end{enumerate}
Since both $\alpha$ and $1-\alpha$ are positive constants, we can thus summarize our observations as follows: The expression $|\phi(t, r)|$ decays as
\[
	| \phi(t, r) | \le 
		\begin{cases}
			\exp(-Q_1(r)) & \text{ for } r \to \infty, \\
			\exp(-Q_2(r)) & \text{ for } r \to - \infty,
		\end{cases}
\]
where the function $Q_1(r) \to \infty$ as $r \to \infty$ at an algebraic rate whereas $Q_2(r) \to \infty$ as $r \to -\infty$ exponentially. Therefore,
the decay of $\phi(t,r)$ to zero as $r \to -\infty$ is much faster than for $r \to \infty$.
\end{rem}

In practical numerical work it will be necessary to approximate the function $\phi(t, \cdot)$ when we actually compute the integral in \eqref{ExDiffusiveIntegral11}. In this context, it is useful to note that we can write 
\[
	\phi(t,r) = \exp(-Q(r)) \Phi(t,r)
\]
where $\Phi(t, \cdot)$ is a bounded function and
\[
	Q(r) = \begin{cases}
			Q_1(r) & \text{ for } r \ge 0, \\
			Q_2(r) & \text{ for } r < 0,
		\end{cases}
\]
with the functions $Q_1$ and $Q_2$ mentioned in Remark \ref{rem:freud-erdos}. From an approximation theory point of view, this is a classical problem; in fact, one would here say that the weight function $\exp(-Q)$ is of Freud type for $r \to \infty$ and of Erd\H os type for $r \to -\infty$ \citep{LL}, so the approximation for $r \in [0,\infty)$ needs to be handled in a different way than for $r \in (-\infty, 0)$.

\subsection{Kernel Function Approximation by Exponential Sums}
Next, we proceed to derive an exponential sum approximation for the kernel function of the integral presented in (\ref{ExDiffusiveIntegral11}). To achieve this, we employ the trapezoidal rule to discretize the integral over an unbounded domain. This rule offers an explicit discretization of the integral over the unbounded domain, expressed as a sum of exponential terms. Importantly, our integrand defined in (\ref{ExDiffusiveIntegralpart2}) exhibits a vital property of exponential decay, as delineated in Theorem \ref{mainthm2}(iv) and discussed in more detail in Remark \ref{rem:freud-erdos}. This characteristic allows for the effective utilization of the well-established trapezoidal rule, known for its accuracy in such scenarios. For a more in-depth understanding of the exponential sum approximation, interested readers can refer to the work presented by \cite{beylkin} and \cite{mclean2018exponential}. 

We recall that it is possible to express the value of $I^{\alpha}_af(t)$ in the form
\begin{equation*}
I^{\alpha}_af(t)=c_{\alpha}\int_a^t\int_{0}^{\infty}\bigg(\frac{u}{t-\tau}\bigg)^{(1-\alpha)}\mathrm e^{-u}u^{-1}f(\tau) \, \mathrm du \, \mathrm d\tau,
\end{equation*}
see, e.g., \citet[proof of Theorem 1]{diethelm2023diffusive}.
By making the substitution $u=(t-\tau)\mathrm e^{r-\mathrm e^{-r}} = (t - \tau) \psi(r)$, we arrive at
\begin{align*}\label{integralpart2}
I^{\alpha}_af(t)= {} &c_{\alpha}\int_a^t\int_{-\infty}^{\infty}(1+\mathrm e^{-r})[\exp(r-\mathrm e^{-r})]^{(1-\alpha)}\nonumber\\
&\quad\times\exp[-(t-\tau)\exp(r-\mathrm e^{-r})]f(\tau) \, \mathrm dr \, \mathrm d\tau
\end{align*}
which can be expressed as
\begin{align*}
I^{\alpha}_af(t)=c_{\alpha}\int_a^tK(t-\tau)f(\tau) \, \mathrm d\tau
\end{align*}
where the kernel $K(t)$ is given by
\begin{eqnarray}\label{Kerneldiff2}
K(t)=&\int_{-\infty}^{\infty}(1+\mathrm e^{-r})[\exp(r-\mathrm e^{-r})]^{(1-\alpha)}\nonumber\\
&\quad\times\exp[-t\exp(r-\mathrm e^{-r})] \, \mathrm dr.
\end{eqnarray}
Applying the trapezoidal rule to discretize the integral with step size $h>0$, we obtain 
\begin{equation}
\begin{array}{ll}
K(t)&\approx \,h\sum_{l=-\infty}^{\infty}(1+\mathrm e^{-lh})[\exp(l h-\mathrm e^{-l h})]^{(1-\alpha)}\nonumber\\
&\quad\quad\quad\quad\times\exp[-t\exp(l h-\mathrm e^{-l h})]\nonumber\\
&=\sum_{l=-\infty}^{\infty} w_l \mathrm e^{-\beta_l t}\label{6}
\end{array}
\end{equation}
with
\begin{align*}
w_l = h(1+\mathrm e^{-lh})[\exp(lh-\mathrm e^{-lh})]^{(1-\alpha)}
\end{align*}
and 
\begin{align*}
\beta_l=\exp(lh-\mathrm e^{-lh}).
\end{align*}

If we limit $t$ to a closed interval $[\delta, b-a]$ with $0 < \delta < b-a < \infty$, we can expand and adapt this method to develop a finite exponential sum approximation for the kernel $K$ given by
\begin{eqnarray}\label{finitesumappr1}
K(t)\approx\sum_{l=-M}^{N} w_l \mathrm e^{-\beta_lt}\quad \text{ for }\,\,t\in [\delta,b-a]
\end{eqnarray}
for suitable choices of $M>0$ and $N>0$.

In \cite{beylkin}, the authors discussed an innovative reduction algorithm based on Prony's method, specifically designed for cases with an excess of terms and small exponents in the exponential sum approximation to a kernel function. This suggests that employing Prony's method allows for a substantial reduction in the number of terms in the exponential sum, as depicted in (\ref{finitesumappr1}), without compromising accuracy. Consequently, numerical algorithms designed using this exponential sum approximation for approximating the fractional integral demonstrate markedly enhanced efficiency when compared to existing numerical methods.

\subsection{Truncation Error Estimation of the Kernel Function}
When dealing with practical computations, especially in numerical analysis and real-world applications, we always need to work with finite sums due to limitations in computational resources. Hence we need to truncate the exponential sum approximation of the kernel function defined in (\ref{6}). In this section, we estimate the two parts of the associated truncation error.

For this, we define an entire analytic function $g(z)$ for $z\in \mathbb{C}$ (the set of complex numbers) and $t>0$ such as
\begin{eqnarray}\label{g}
g(z)=(1+\mathrm e^{-z})\exp[(1-\alpha)(z-\mathrm e^{-z})]\nonumber\\
\,\times\exp[-t\exp(z-\mathrm e^{-z})]  
\end{eqnarray}
which is the analytical extension of the integrand function defined in equation (\ref{Kerneldiff2}). Then, we have the following result which estimates the bounds of $g(z)$:

\begin{thm}
\label{thm:truncation}
Let $h>0$, $0<\theta<1$ and $t\in (0,1]$. Then the function $g$ defined in (\ref{g}) satisfies the bounds
\begin{equation}\label{7}
h\sum_{l=-\infty}^{-M-1}g(lh)\leq C_1 \exp((\alpha-1)\mathrm e^{Mh})
\end{equation}
for
\[
Mh \geq
               \begin{cases}
                  \log(\alpha(1-\alpha)^{-1}), & \text{ for } 1/2<\alpha<1, \\
                   0, & \text{ for } 0 <\alpha \leq 1/2.
                \end{cases}
\]
and 
\begin{equation}\label{8}
h\sum_{l=N+1}^{\infty}g(lh)\leq \frac{C_2}{t^\beta} \exp(-\theta t\mathrm e^{Nh-1}) 
\end{equation}
for
\begin{align*}
Nh\geq 1+\log((1-\alpha)t^{-1}).
\end{align*}
For $0\leq\alpha<1$, (\ref{8}) also holds for $\theta=1$.
\end{thm}

\begin{pf}
The proof can be directly derived from Lemma 2 in \cite{mclean2018exponential} with appropriate modifications.
\qed
\end{pf}

Note that formally the statement of Theorem \ref{thm:truncation} is limited to $t \le 1$. However, this limitation can be overcome by a simple rescaling technique as indicated in \citet[p.~928]{mclean2018exponential}.

We can readily verify that the number $\Lambda=M+1+N$ of terms needed to guarantee a relative accuracy $\epsilon$ for $\delta\leq t\leq 1$ is of order $(\log \epsilon^{-1})\log(\delta^{-1}\log\epsilon^{-1})$.

\section{Efficient Numerical Algorithms for Approximate Evaluation of Fractional Integral}\label{Numerical Algorithms}
In this section, we establish two distinct efficient numerical algorithms for the approximate evaluation of the fractional integral represented in (\ref{ExDiffusiveIntegral11}). The first approach involves utilizing an exponential sum approximation technique on the kernel function, as defined in (\ref{finitesumappr1}). The second approach employs the Gauss-Laguerre quadrature formula.

\subsection{Numerical Algorithm for Approximating Fractional Integral using Exponential Sum Approximation of the Kernel Function}
For the numerical computation of $I_a^{\alpha}f(t)$ at the points of the grid $a=t_0<t_1<\cdots <t_{P}=b$, let $F^n = f(t_n)$ and define a piecewise constant interpolant $\tilde{F}(t)=F^n$ for $t\in J_n=(t_{n-1},t_n)$. Then
\begin{equation}
  I^{\alpha}_af(t_n) =c_{\alpha}\int_a^{t_n}K(t_n-\tau)f(\tau) \, \mathrm d\tau \approx \sum_{j=1}^{n}z_{nj}F^j
\end{equation}
where 
\begin{equation}
z_{nj}=c_{\alpha}\int_{J_j}K(t_n-\tau) \, \mathrm d\tau.
\end{equation}

The computational complexity for calculating this sum within the range of $1\leq n\leq P$ in a conventional way scales proportionally to $P^2$. Such quadratic growth can pose challenges, particularly in scenarios where each $F^j$ represents a large vector rather than a scalar. Furthermore, storing $F^j$ in active memory for all time levels $j$ may not be feasible. 

Here we define an efficient algorithm to avoid these challenges based on the exponential sum approximation to the kernel $K$ such as 
\begin{equation}
K(t)\approx \sum_{l=1}^{\Lambda}w_l \mathrm e^{\beta_l t} \quad \textit{for}\quad \delta\leq t\leq b-a, 
\end{equation} 
provided that the moderate number of terms $\Lambda$ can achieve adequate accuracy for a choice of $\delta>0$ that is smaller than the time step $\Delta t_n=t_n-t_{n-1}$ for all $n$. Certainly, if $\Delta t_n\geq \delta$ then $\delta\leq t_n-\tau \leq b-a$ for $0\leq \tau \leq t_{n-1}$, and therefore
\begin{eqnarray*}
\sum_{j=1}^{n-1}z_{nj}F^j &=&c_{\alpha}\int_{0}^{t_{n-1}}K(t_n-\tau)\tilde{F}(\tau) \, \mathrm d\tau\\
&\approx& c_{\alpha}\int_{0}^{t_{n-1}}\sum_{l=1}^{\Lambda}w_l \mathrm e^{\beta_l (t_n-\tau)}\tilde{F}(\tau) \, \mathrm d\tau\\
&=& \sum_{l=1}^{\Lambda} \Phi_l^n
\end{eqnarray*}
where 
\begin{equation*}
\Phi_l^n  =  c_{\alpha}w_l \int_{0}^{t_{n-1}}\mathrm e^{\beta_l (t_n-\tau)}\tilde{F}(\tau) \, \mathrm d\tau
 	=   \sum_{j=1}^{n-1}K_{lnj}F^j
\end{equation*}
and 
\begin{equation}
K_{lnj}= c_{\alpha}w_l\int_{J_j}\mathrm e^{\beta_l (t_n-\tau)} \, \mathrm d\tau.
\end{equation}
Hence the approximated value of the integral is expressed as
\begin{equation}\label{Trapappro1}
I^{\alpha}_a\tilde{F}(t_n)\approx z_{nn}F^n+ \sum_{l=1}^{\Lambda} \Phi_l^n
\end{equation}
and evaluated utilizing the recursion formula
\begin{equation}\label{Trapappro2}
\Phi_l^n=K_{l,n,n-1}F^{n-1}+\mathrm e^{\beta l \Delta t_n}\Phi_l^{n-1} \qquad (2\leq n \leq P)
\end{equation}
with starting value $\Phi_l^1=0$.

Thus we can efficiently compute $I^{\alpha}_a\tilde{F}(t_n)$ with a satisfactory level of accuracy using a number of operations that is proportional to $\Lambda \cdot P$. This approach yields substantial computational savings, particularly when $L \ll P$. Moreover, we have the flexibility to overwrite $\Phi_l^{n-1}$ with $\Phi_l^{n}$, and $F^{n-1}$ with $F^{n}$, effectively reducing the active storage requirement from being proportional to $P$ to being proportional to $\Lambda$.

\subsection{Numerical Algorithm for Approximating Fractional Integral using Gauss-Laguerre Quadrature Formula}\label{GLQuadrature}

Next, we proceed to approximate the integral (\ref{ExDiffusiveIntegral11}) using the $\Lambda$-point Gauss-Laguerre quadrature formula. This particular quadrature method is specifically designed for evaluating singular integrals, particularly those with the weight function $\mathrm e^{-u}$ over the interval $[0,\infty)$.

We decompose the integral as follows:
\begin{eqnarray*}
I^{\alpha}_af(t) &=& \int_{-\infty}^{\infty}\phi(t,r) \, \mathrm dr\\
&= &\int_{-\infty}^{0}\phi(t,r)dr +\int_{0}^{\infty}\phi(t,r) \, \mathrm dr\\
&=&J_1+J_2.
\end{eqnarray*}
Substituting $r=\frac{-s}{1-\alpha}$ in $J_1$ and $r=\frac{s}{\alpha}$ in $J_2$, we obtain the following expression:
\begin{equation*}
\begin{split}
 I^{\alpha}_af(t) =& \frac{1}{1-\alpha}\int_{0}^{\infty}\mathrm e^{-s}\mathrm e^{s}\phi(t,-s/(1-\alpha)) \, \mathrm ds\\
&+\frac{1}{\alpha}\int_{0}^{\infty}\mathrm e^{-s}\mathrm e^{s}\phi(t,s/\alpha) \, \mathrm ds.
\end{split}
\end{equation*}
By introducing the function $\hat{\phi}(t,s)$ defined as:
\begin{eqnarray*}
\hat{\phi}(t,s) =& \mathrm e^{s} \bigg(\frac{1}{1-\alpha}\phi(t,-s/(1-\alpha))+\frac{1}{\alpha}\phi(t,s/\alpha)\bigg),
\end{eqnarray*}
We can then express the integral as:
\begin{eqnarray*}
I^{\alpha}_af(t) =\int_{0}^{\infty}\mathrm e^{-s}\hat{\phi}(t,s) \, \mathrm ds \approx  Q_\Lambda^{GLa}[\hat{\phi}(t,\cdot)],
\end{eqnarray*}
where the expression
\begin{eqnarray*}
Q_\Lambda^{GLa}[\hat{\phi}]=\sum_{l=1}^{\Lambda}w_l^{GLa}\hat{\phi}(x_l^{GLa})
\end{eqnarray*}
represents the $\Lambda$-point Gauss-Laguerre quadrature formula with weights denoted as $w_l^{GLa}$ and nodes as $x_l^{GLa}$.

To determine the Gauss-Laguerre nodes $x_l^{GLa}$ and their corresponding weights $w_l^{GLa}$, we follow these steps:
\begin{itemize}
  \item Solve $L_\Lambda(x_l^{GLa})=0$ for $n=1,2,\ldots, \Lambda$ to obtain the Gauss-Laguerre nodes $x_l^{GLa}$.
  \item To calculate the weights, use the formula 
\begin{eqnarray*}
w_l^{GLa}=\frac{x_l^{GLa}}{[L_{\Lambda+1}(x_l^{GLa})^2]},
\end{eqnarray*}
\end{itemize}
where $L_{\Lambda+1}$ represents the Laguerre polynomial of order $\Lambda+1$.

For a more detailed understanding of the Gauss-Laguerre quadrature formula, we refer to \cite{davis2007methods}.

\begin{thm}
Under the assumptions specified in Theorem (\ref{mainthm2}), we establish the following convergence result:
\begin{equation*}
\lim_{\Lambda \rightarrow \infty} Q_\Lambda^{GLa}[\hat{\phi}(t,\cdot)]=I^{\alpha}_af(t)
\end{equation*}
for all $t\in [a,b]$.
\end{thm}

\begin{pf}
Leveraging the properties presented in Theorem \ref{mainthm2}(iii) concerning the multiple differentiability of the function $\phi(t,\cdot)$, and coupled with the decay characteristics outlined in Theorem \ref{mainthm2}(iv), we are equipped to apply a well-established convergence result associated with the Gauss-Laguerre quadrature formula, as documented in \cite{davis2007methods}. Consequently, we derive the desired result.
\qed
\end{pf}

We emphasize here that, in view of the asymptotic behaviour of the integrand established in Remark \ref{rem:freud-erdos}, the choice of the Gauss-Laguerre method is natural for the integral $J_2$. i.e.\ for the contribution of the integral (\ref{ExDiffusiveIntegral11}) that stems from the subinterval $[0, \infty)$. For the remaining component $J_1$, we had seen a much faster decay to zero in Remark \ref{rem:freud-erdos}. Therefore, one might either follow the strategy suggested by \cite{MM} and truncate the quadrature sum, i.e.\ use only the summands $\ell = 1, 2, \ldots, \Lambda^*$ with $\Lambda^* \ll \Lambda$ (a typical choice might be $\Lambda^* = \Lambda/2$) to reduce the computational effort and simultaneously improve the accuracy, or use a Gauss-type quadrature formula for the Erd\H os-type weight function \citep{LL} instead of the Laguerre weight.

\subsection{Methodology and Computational Procedures}\label{GLQuadratureimplimentation}
We are now prepared to outline the method we propose for numerically computing $I^{\alpha}_a f(t_n)$, where $n = 1, 2, \ldots, P$. In this algorithm, we employ the symbol $\hat{\phi}_l$ to denote the approximate value of $\hat{\phi}(x_l^{GLa},t_n)$ for the current time step, which corresponds to the currently evaluated $n$ value.

Given the initial point $a\in \mathbb{R}$, the order $\alpha\in (0,1)$, the grid points $t_n$ for $n = 1, 2, \ldots, P$, and the number of quadrature nodes $\Lambda \in \mathbb{N}$,

\begin{enumerate}
  \item For $l= 1, 2, \cdots, \Lambda$:
  \begin{description}
    \item[a.] Compute the Gauss Laguerre nodes $x_l^{GLa}$ and the associated weights $w_l^{GLa}$.
    \item[b.] Define the auxiliary terms $r_l\leftarrow \frac{-x_l^{GLa}}{1-\alpha}$ and $\tilde{r}_l\leftarrow \frac{x_l^{GLa}}{\alpha}$
    \item[c.] Initialize $\phi_l\leftarrow 0$ and  $\tilde{\phi}_l\leftarrow 0$ to denote the initial condition for the differential equation (\ref{DE})
  \end{description}
  \item For $n = 1, 2, \ldots, P$: 
   \begin{description}
    \item[a.] set $h\leftarrow t_n-t_{n-1}$.
    \pagebreak
    \item[b.] For $l= 1, 2, \ldots, \Lambda$:
     \begin{description}
    \item[(i)] Update the value of $\phi_l$ by solving the corresponding differential equation (\ref{DE2}) using the backward Euler method
   \begin{eqnarray}\label{4}
    \phi_l\leftarrow &\frac{1}{1+h\mathrm e^{(r_l-\mathrm e^{-r_l})}}\bigg[\phi_l+hc_{\alpha}(1+\mathrm e^{-r_l})\nonumber\\
&\times \mathrm e^{(1-\alpha)(r_l-\mathrm e^{-r_l})}f(t_n)\bigg]
    \end{eqnarray}
    \item[(ii)] Similarly, update the value of $\tilde{\phi}_l$ by
       \begin{eqnarray}\label{5}
    \tilde{\phi}_l\leftarrow &\frac{1}{1+h\mathrm e^{(\tilde{r}_l-\mathrm e^{-\tilde{r}_l})}}\bigg[\tilde{\phi}_l+hc_{\alpha}(1+\mathrm e^{-\tilde{r}_l})\nonumber\\
&\times \mathrm e^{(1-\alpha)(\tilde{r}_l-\mathrm e^{-\tilde{r}_l})}f(t_n)\bigg]
    \end{eqnarray}
  \end{description}
  \item[c.]  Calculate the desired approximate value for $I^{\alpha}_a f(t_n)$ using the formula
     \begin{equation*}
   I^{\alpha}_a f(t_n)=\sum_{l=1}^{\Lambda}w_l^{GLa}\exp(x_l^{GLa})\bigg[\frac{1}{1-\alpha}\phi_l+\frac{1}{\alpha}\tilde{\phi}_l\bigg].
    \end{equation*}
  \end{description}
\end{enumerate}

In this context, we opt for the backward Euler method, taking into account the constant factor that multiplies the unknown function $\phi(\cdot,r)$ in equation (\ref{DE2}). This choice is motivated by the need to utilize an A-stable method, and the backward Euler method represents the simplest among such options. Alternatively, one could consider the trapezoidal method, which is also an A-stable approach. However, implementing the trapezoidal method would require modifications to the formulas provided in (\ref{4}) and (\ref{5}), as follows:

\begin{eqnarray}\label{algo5}
    \phi_l\leftarrow \frac{1}{1+\frac{h}{2}\mathrm e^{(r_l-\mathrm e^{-r_l})}}\bigg[\bigg(1-\frac{h}{2}\mathrm e^{(r_l-\mathrm e^{-r_l})}\bigg)\phi_l+\frac{h}{2}\nonumber\\
    \times c_{\alpha}(1+\mathrm e^{-r_l})\mathrm e^{(1-\alpha)(r_l-\mathrm e^{-r_l})}[f(t_n)+f(t_{n-1})]\bigg]
    \end{eqnarray}
and
\begin{eqnarray}\label{algo6}
    \tilde{\phi}_l\leftarrow \frac{1}{1+\frac{h}{2}\mathrm e^{(\tilde{r}_l-\mathrm e^{-\tilde{r}_l})}}\bigg[\bigg(1-\frac{h}{2}\mathrm e^{(\tilde{r}_l-\mathrm e^{-\tilde{r}_l})}\bigg)\phi_l+\frac{h}{2}\nonumber\\
    \times c_{\alpha}(1+\mathrm e^{-\tilde{r}_l})\mathrm e^{(1-\alpha)(\tilde{r}_l-\mathrm e^{-\tilde{r}_l})}[f(t_n)+f(t_{n-1})]\bigg]
    \end{eqnarray}
    respectively.

\section{Conclusion and Future Directions}

In this research, we have introduced an innovative diffusive representation method for the fractional integral. This novel approach offers a practical and efficient means of numerically approximating the integral. By strategically selecting a specific form of the admissible transformation function, we've leveraged the advantageous property of exponential decay within our integrand. This feature has allowed us to effectively employ the trapezoidal rule for discretization, even over unbounded domains, or the widely-used Gauss-Laguerre quadrature method. Importantly, this approach harmonizes with common techniques applied in solving problems featuring integer-order operators, yielding benefits in memory efficiency and computational simplicity.

While this study marks a significant advancement in the field of numerical approximation for fractional integrals, there are promising directions for future research and development. These include broadening the scope of the diffusive representation to encompass a wider range of admissible transformation functions and domains, thereby enhancing its versatility. Additionally, ongoing efforts can focus on optimizing the numerical algorithms, particularly in scenarios involving intricate functions. The practical application of the diffusive representation in various real-world domains, such as environmental modeling, finance, and control systems, will help gauge its efficacy. It is imperative to conduct rigorous error analysis to validate the accuracy and convergence properties of the proposed numerical methods for diverse applications.


\bibliography{ifacconf}             
                                                   







\end{document}